\documentclass[11pt]{article}
\usepackage{amsmath,amssymb,epsf}
\usepackage[french]{babel}
\usepackage[applemac]{inputenc}
\advance \textwidth 3cm
\advance \textheight 3cm
\usepackage{amsfonts}
\usepackage{latexsym}
\newtheorem{theoreme}{Th\'eor\`eme}
\newtheorem{lemme}{Lemme}
\newtheorem{corollaire}{Corollaire}

\newtheorem{proposition}{Proposition}
\newtheorem{prop}{Proprit}

\newenvironment{Preuve}[1]{\par\noindent\underline{Preuve #1} :\quad}%
{\unskip\nobreak\hfil\penalty50\hskip2em\null\nobreak\hfil%
$\Box$\parfillskip0pt\par\medskip}%versions F et A

\newcommand{\trace}{\mathrm{Tr}}

\title{Expression asymptotique des valeurs propres d'une matrice de Toeplitz \`a symbole r\'eel.}

\author{ Philippe Rambour\thanks{Universit\'{e} de Paris Sud,
      B\^atiment 425; F-91405
Orsay Cedex;
tel : 01 69 15 57 28 ; fax 01 69 15 60 19
      \mbox{e-mail : philippe.rambour@math.u-psud.fr}}
      }
\date{}
\begin{document}
\maketitle
%\noindent
  \renewcommand{\abstractname}{R\'ESUM\'E}
   \begin{abstract}
    \textbf{Expression asymptotique des valeurs propres d'une matrice de Toeplitz \`a symbole r\'eel.}\\
  Ce travail donne deux r\'esultats que nous obtenons \`a partir 
  d'une formule d'inversion des matrices de Toeplitz \`a symbole r\'eel, qui a \'et\'e \'etablie dans un pr\'ec\'edent article.
 Le premier de ces r\'esultats donne une expression asymptotique des valeurs propres minimales d'une matrice de Toeplitz 
$T_{N}\left( f\right)$ o\`u $f$ est une fonction  p\'eriodique, paire et d\'erivable sur $[0, 2 \pi[$. 
Ensuite, nous d\'emontrons qu'une matrice de Toeplitz bande de taille $N\times N$ dont le symbole n'admet pas de z\'eros 
complexes de module un est inversible et que son inverse quand $N$ tend vers l'infini est une matrice bande, \`a des termesn\'egligeables pr\`es. Une cons\'equence de ce dernier r\'esultat  est une estimation asymptotique des coefficients de l'inverse d'une matrice de Toeplitz \`a symbole r\'egulier. 
\end{abstract}
 \renewcommand{\abstractname}{ABSTRACT}
                 \begin{abstract}
                   \textbf{Asymptotic expression of the eigenvalues of a Toeplitz matrix with a real symbol.}\\
                This work provides two results obtained as a consequence of an inversion formula for 
                Toeplitz matrices with real symbol. First we obtain an asymptotic expression for  the minimal eigenvalues of a Toeplitz 
                matrix with a symbol which is periodic, even and derivable on $[0,2 \pi[$. Next we prove that a Toeplitz band matrix with a symbol without zeros on the united circle is invertible with an inverse which is essentially a band matrix. As 
                a consequence of this last statement we give an asymptotic esimation for the entries of the inverse of a Toeplitz matrix with a regular symbol.\\
 \end{abstract}

%\vskip 30 mm
\textbf{\large{Mots clef}}\\
\textbf{Matrices de Toeplitz, valeurs propres , formule d'inversion, matrices de Toeplitz bandes.}
\section {Introduction}
En toute g\'en\'eralit\'e, on dit qu'une matrice $N\times N$ est une matrice
  de Toeplitz d'ordre  $N$ s'il existe $2N-1$ nombres complexes 
  $a_{-(N-1)}, a_{-(N-2)}, \cdots, a_{-1}, a_{0}, a_{1}, \cdots, a_{N-2},a_{N-1}$ 
  tels que $\left(T_{N}\right)_{i,j}= a_{j-i}$, pour $1\le i,j \le N$.\\
 Rappelons  que si $f$ est une fonction de $L^1(\mathbb T)$ (avec 
 $\mathbb T = \mathbb R / 2\pi \mathbb Z$), 
 on appelle matrice de Toeplitz d'ordre $N$ de symbole $f$, 
 et on note $T_N(f)$, la matrice $(N+1) \times (N+1)$ telle que 
 $\left(T_N (f)\right) _{k+1,l+1} =\hat f (k-l) \quad \forall \,k,l \quad 
 0\le k,l \le N$ o  $\hat h (j)$ dsigne le coefficient de Fourier 
 d'ordre $j$ d'une fonction $h$ (une bonne r\'ef\'erence peut \^{e}tre \cite{Bo.3}). Dans 
 la premi\`ere partie de ce travail, on s'int\'eresse   l'expression asymptotique lorsque $N$ tend vers l'infini  des valeurs propres des matrices de Toeplitz $T_{N}(f)$ o\`u
 $ f$ est une fonction p\'eriodique, paire, d\'erivable. 
 \\
Quand on a affaire \`a un symbole positif, born\'e et int\'egrable, il est connu que ces valeurs propres sont \'equidistribu\'ees au sens de Weyl avec les quantit\'es 
 $f(-\pi+\frac{2j\pi}{N+1})$, $1\le j\le N$. Rappelons que deux familles de sous-ensembles de $\mathbb R$,
 $\left(\mu_{j}^{(N)}\right)_{j=1}^N$ et $\left(\kappa_{j}^{(N)}\right)_{j=1}^N$ sont dites \'egalement distribu\'ees au sens de Weyl si pour toute fonction $h$ appartenant 
 \`a $\mathcal C(\mathbb R)$ on a 
 $ \displaystyle{\lim_{N\to\infty}\frac{1}{N}
 \sum_{j=1}^N \left( f(\mu_{j}^{(N)}) -f(\kappa_{j}^{(N)})\right)=0}$ (on peut consulter 
 \cite{Bo.3} et \cite{GS} sur ces questions). 
 Notre article compl\`ete des travaux pr\'ec\'edents (\cite {RS020}) o\`u l'on a consid\'er\'e  des 
 fonctions $f\in L^1(\mathbb T)$ paires, positives et strictement monotones sur $[0, \pi]$. Le th\'eor\`eme \ref{principal} est une extension de ce travail \`a des fonctions qui sont  paires et d\'erivables, mais plus n\'ecessairement monotones. Ceci peut \^{e}tre reli\'e \`a ce qu'on sait sur les valeurs propres de $T_{N}(\vert P\vert^2)$ lorsque $P$ est un polyn\^{o}me trigonom\'etrique
de degr\'e un. Ces valeurs propres sont alors de la forme 
$\vert P( \xi + \frac{j \pi}{N+1})\vert^2$, $1\le j\le N+1$  (\cite{GS}) o\`u $\vert P(\xi)\vert^2$ est le minimum de $\vert P\vert^2$. Le th\'eor\`eme \ref{principal} est \'egalement conforme \`a la conjecture de Whittle (\cite{WHI},\cite{BRDa}), qui annonce que pour 
  certains des symboles positifs $ \theta\mapsto h(e^{i\theta})$ v\'erifiant les hypoth\`eses du th\'eor\`eme \ref{principal}  par exemple les densit\'es spectrales des processus ARIMA (voir \cite{B}) les valeurs propres de $T_{N}(h)$ sont en fait de la forme $h( e^{i k \pi/N})$.\\
 On sait 
 (\cite{GS}) que si 
 $f\in L^1(\mathbb T)$ est une fonction positive, alors les valeurs propres de $T_{N}(f)$
 sont comprises entre le minimum et le maximum  de $f$ sur le tore. En particulier, la valeur propre minimale de $T_{N}(f)$ tend vers le minimum de $f$ quand
  $N \rightarrow + \infty$. Obtenir dans ce cas une approximation plus pr\'ecise de cette valeur propre minimale est un probl\`eme sp\'ecifique qui a fait l'objet de nombreux travaux. On peut voir par exemple \cite{JMR01},\cite{Part}, \cite{Part2}, \cite{ W3},\cite{SpSt} pour une approche directe, et \cite{Bow}, \cite{RS1111}, \cite{RQuebl} pour une m\'ethode qui consiste \`a calculer la plus grande valeur propre de la matrice inverse.\\
Dans ce travail, nos d\'emonstrations sont  bas\'ees sur une formule d'inversion des matrices de Toeplitz \`a symbole r\'eel, formule qui
  a \'et\'e \'etablie dans \cite{RRQuebl}, qui g\'en\'eralise les m\'ethodes d'inversion classiques concernant les matrices de Toeplitz \`a symbole positif, et qui est mieux 
  adapt\'ee \`a la situation. Une telle g\'en\'eralisation a \'egalement \'et\'e utilis\'ee dans 
 \cite{De.Ri} \\
 Dans la seconde partie de l'article, nous utilisons cette formule pour obtenir un r\'esultat sur l'inverse des matrices de Toeplitz bande 
 $T_{N}$ dont le symbole est une fonction  
 $\theta\mapsto \displaystyle{ \sum_{j=-n_{0}}^{n_{0}}  a_{j} e^{i j\theta}}$ telle que 
 $ \theta\mapsto  \displaystyle{ \sum_{j=-n_{0}}^{n_{0}}  a_{j} e^{i (j+n_{0}) \theta}}$ ne  s'annulant pas en z\'ero (voir le th\'eor\`eme \ref{theoreme4}). Si le nombre de racines de module strictement 
  sup\'erieur \`a $1$ est \'egal au nombre de racines de module strictement 
  inf\'erieur \`a $1$, nous d\'emontrons que l'inverse d'une telle matrice de Toeplitz est asymptotiquement 
  une matrice bande, au sens o\`u lorsque $\vert k-l\vert $ est suffisamment grand, le coefficient 
  $(T_{N})^{-1}_{k,l}$ est d'ordre $\rho^{\vert k-l\vert}$ pour une constante $\rho$ 
  comprise strictement entre $0$ et $1$. Nous obtenons 
  \'egalement une approximation des termes de l'inverse proches de la diagonale.
  Rappelons que de nombreux algorithmes de calcul permettent de
  calculer num\'eriquement l'inverse des matrices de Toeplitz bande, notamment pour des expos\'es r\'ecents, (\cite{Hen-Bar}, \cite{Mal-Sad2}). Comme cons\'equence de ce dernier 
  r\'esultat nous obtenons, gr\^{a}ce \`a une approximation polynomiale, une expression
  asymptotique lorsque $k-l$ tend vers l'infini des coefficients
   $T^{-1}_{N}(\psi)$ o\`u $\psi$ est un symbole r\'egulier, (on appelle fonction r\'eguli\`ere une fonction de  
 $L^1(\mathbb T)$ strictement positive sur le tore $ \mathbb T$). Enfin 
 si $h$ est une fonction positive dans $L^1(\mathbb T)$, 
 le polyn\^{o}me pr\'edicteur de degr\'e $M$ associ\'e \`a $h$ est 
 le polyn\^{o}me trigonom\'etrique 
 $P_{M} : \theta\mapsto \displaystyle{ \sum_{u=0}^M 
 \beta_{u,M} e^{iu\theta}}$ o\`u 
 $\beta_{u,M}= \frac{\left(T_{M}^{-1}(h)\right))_{(u,1)}}
{ \sqrt{ \left(T_{M}^{-1}(h)\right))_{(1,1)}}}. $ 
 On sait que $P_{M}$ v\'erifie la propri\'et\'e suivante :
 \begin{prop}
 pour tout entier $s$, $-M\le s\le M$ on a 
 $$ \hat h (s) = \widehat{\frac{1}{\vert P_{M}\vert ^2}} (s).$$
 \end{prop}
 Pour les polyn\^{o}mes pr\'edicteurs, le lecteur pourra se r\'ef\'erer \`a  \cite{Ld}.
 \\
    Dans le paragraphe suivant, nous allons rappeler les notations utilis\'ees 
   et expliciter la formule d'inversion utilis\'ee ici.
  \section{Formule d'inversion}
  \subsection{Rappel des notations}
 Dans toute la suite, nous noterons par $\chi$ la fonction de $[0,2 \pi[$ dans $\mathbb C$ 
 d\'efinie par $\theta \mapsto e^{i\theta}$ et nous noterons $\mathcal P_{N}$ l'ensemble 
 des polyn\^{o}mes trigonom\'etriques de degr\'e inf\'erieur ou \'egal \`a $N$ et $\mathcal P_{[-N,N]}$ le sous-espace vectoriel de $L^{1}( \mathbb T$ engendr\'e par les fonctions $\chi^{j} -N\le j\le N$. Nous consid\'ererons \'egalement les sous-espaces suivants de $L^2(\mathbb T)$ (espaces de Hardy),
  \begin{enumerate}
  \item
  $$ H^+ (\mathbb T) = \{\psi \in \mathbb T \vert u<0 \Rightarrow \hat{\psi} (u)=0\},
  $$
  \item
  $$
 \left(  H^+ (\mathbb T)\right)^\bot = \{\psi \in \mathbb T \vert u\ge 0 \Rightarrow \hat{\psi} (u)=0\},$$
 \item
 $$
 H^{\infty}(\mathbb T) = H^+ (\mathbb T) \cap L^\infty (\mathbb T).
  $$
  o\`u $ L^\infty (\mathbb T)$ est l'espace des fonctions mesurables born\'ees presque partout pour la mesure de Lebesgue, muni de la distance $d_{\infty}$ associ\'ee \`a la norme $\Vert . \Vert_{\infty}$ 
  d\'efinie pour toute fonction $\psi$ dans $L^{1}(\mathbb T)$ par $\Vert \psi \Vert_{\infty} =\sup_{\theta\in [0, 2\pi[} \vert f(\theta)\vert$.
  \end{enumerate}
 On peut alors d\'efinir $\pi_{N},\pi_{+},\pi_{-},$ les projections respectives de 
 $L^2 (\mathbb T) $ sur $\mathcal P_{N}$, $H^+ (\mathbb T)$ et 
  $\left(H^+(\mathbb T)\right)^\bot $. \\
  Si maintenant $g_{1}$ est un \'el\'ement de $H^+ (\mathbb T)$ et $g_{2}$ un 
  \'el\'ement tel que  $\bar g_{2} \in H^+ (\mathbb T)$, on pose $\Phi_{N}= \chi^{N+1} \frac{g_{1}}{g_{2}}$
  et $\tilde \Phi_{N}= \chi^{-N-1} \frac{g_{2}}{g_{1}}$ et on d\'efinit les op\'erateurs de Haenkel 
 $H_{\Phi_{N}}$ et $H_{\tilde \Phi_{N}}$ par 
 \begin{align*} 
 H_{\Phi_{N}} : \, & H^+ (\mathbb T) \rightarrow \left(H^+(\mathbb T)\right)^\bot \\
                          & \psi \mapsto \pi_{-} (\Phi_{N}\psi) 
 \end{align*}
  et 
   \begin{align*} 
 H_{\tilde \Phi_{N}} : \, & \left(H^+(\mathbb T)\right)^\bot \rightarrow H^+ (\mathbb T) \\
                                  & \psi \mapsto \pi_{+} (\tilde \Phi_{N}\psi). 
 \end{align*}
    \subsection{La formule d'inversion proprement dite}
  Classiquement, on sait calculer explicitement les coefficients de l'inverse  d'une matrice de Toeplitz de symbole positif $f$ qui v\'erifie $\ln f \in L^1 (\mathbb T)$. Ces formules 
  utilisent fortement le fait qu'une telle fonction $f$ peut s'\'ecrire $f = g \bar g$ avec 
  $g \in H^+ (\mathbb T)$  (\cite{RRS}). Certains probl\`emes, notamment la recherche des valeurs
  propres, n\'ecessitent d'inverser des matrices de Toeplitz dont le symbole 
  n'est pas n\'ecessairement positif (voir, par exemple, \cite{JMR01} \cite {RS020}).
  Pour ce faire, nous allons consid\'erer des fonctions $f$ qui ne s'annulent pas sur le tore  admettant une d\'ecomposition  $f=g_{1}g_{2}$, avec 
   $g_{1}, g_{1}^{-1} \in H^+(\mathbb T)$, $\bar g_{2}, {\bar g_{2}}^{-1}\in H^+(\mathbb T)$. 
  Nous pouvons alors \'enoncer le th\'eor\`eme :
  \begin{theoreme}\label{UN}
  Si la fonction $f$ est comme ci-dessus, on suppose que les fonctions $g_{1}$ et $g_{2}$ v\'erifient,  
  $$ 
  \lim_{N\rightarrow + \infty} d_{\infty}
  \left( \frac{1}{\vert g_{2}\vert ^{2}}, \mathcal P_{[-N,N}\right) =0
  \quad 
  \mathrm{et}
  \quad 
  d_{\infty}
  \left( \frac{1}{\vert g_{1}\vert ^{2}}, \mathcal P_{[-N,N]}\right) < + \infty
  $$
  ou
   $$ 
  \lim_{N\rightarrow + \infty} d_{\infty}
  \left( \frac{1}{\vert g_{1}\vert ^{2}}, \mathcal P_{[-N,N]}\right) =0
  \quad 
  \mathrm{et}
  \quad 
  d_{\infty}
  \left( \frac{1}{\vert g_{2}\vert ^{2}}, \mathcal P_{[-N,N]}\right) < \infty
$$
  alors, pour $N$ suffisamment grand,
  \begin{enumerate}
  \item
   l'in\'egalit\'e $\Vert H_{\tilde \Phi_{N}} H_{\Phi_{N}} \Vert_{2}<1$ est v\'erifi\'ee (ce qui assure bien s\^{u}r l'inversibilit\'e de l'op\'erateur 
  $I - H_{\tilde \Phi_{N}} H_{\Phi_{N}} $),
  \item
  la matrice $T_{N}(f)$ est inversible,
  \item
    pour tout polyn\^{o}me $Q$ dans $\mathcal P_{N}$, on a 
 $$
  T_{N}^{-1} (f) (Q) =
  \pi_{+}(Q g_{2}^{-1}) g_{1}^{-1} - \pi_{+}
  \left ( \left(( I- H_{\tilde \Phi_{N}} H_{\Phi_{N}} )^{-1} \pi_{+} 
  \left( \pi_{+} (Q g_{2}^{-1}) \tilde \Phi_{N}\right)\right) \Phi_{N}\right ) g_{1}^{-1}.$$
  \end{enumerate}
  \end{theoreme} 
  Ce th\'eor\`eme a \'et\'e obtenu dans \cite{ RRQuebl} pour inverser des matrices \`a blocs que nous interpr\'etons comme des matrices de Toeplitz tronqu\'ees dont le symbole 
  est une fonction matricielle contenue dans $L^2_{\mathcal M_{n}} (\mathbb T)
  =\{ M : \mathbb T \rightarrow \mathcal M_{n} (\mathbb C) \Big\vert 
 \int_{\mathbb T}\Vert M\Vert_{2}^2 d\sigma< + \infty\},$ o\`u 
 $\Vert M \Vert_{2}= [ \trace (M M ^*]^{1/2}$ et o\`u $\sigma$ est la mesure de Lebesgue su le tore $\mathbb T$. N\'eanmoins, il est facile de v\'erifier 
  qu'avec les hypoth\`eses que nous nous sommes donn\'ees (voir \cite{ RRQuebl}, lemme 7.1.5 et corollaire 7.1.6), la d\'emonstration se transpose
  sans difficult\'e au cas des matrices de Toeplitz dont le symbole est une fonction de 
  $\mathbb T$ dans $\mathbb C$. Nous allons appliquer cette formule 
  pour inverser des matrices de Toeplitz d'ordre $N$ dont le symbole $f$ v\'erifie 
  $f (\chi) =  Cg_{1}(\chi)g_{2}(\bar \chi)$  o\`u $C$ est une constante et 
  $g_{1}$ et $g_{2}$ deux fractions rationnelles \`a coefficients complexes dont les z\'eros et les p\^{o}les sont de module strictement inf\'erieur \`a $1$. 
  Nous avons alors facilement la proposition suivante, que nous utiliserons pour \'etablir nos r\'esultats.
 \begin{proposition}   Les fonctions $g_{1}$ et $g_{2}$ \'etant dans $C^{\infty} (\mathbb T)$ elles v\'erifient 
    les hypoth\`eses du th\'eor\`eme \ref{UN}, et donc $T_{N}(f$) est 
    inversible.
    \end{proposition}
       \section{Une application aux valeurs propres des matrices de Toeplitz}
       \begin{theoreme}\label{principal}
    On consid\`ere une fonction $f$, paire, $2\pi$ p\'eriodique, de classe $C^\infty\left( [0, \pi]\right)$.
 Pour tout $\lambda$ dans $ [\min_{\theta \in [0,\pi]}, \max_{\theta \in [0,\pi]}]$ on suppose qu'il existe 
 $r_{\lambda}$ r\'eels, distincts ou non, not\'es $\theta_{1},
 \ldots, \theta_{r_{\lambda}}$ et 
 une fonction $G_{\lambda}$ continue et strictement positive sur $[-\pi, \pi]$ tel que, pour tout $\theta \in [0, \pi]$ on ait 
 $$ f(\theta-\lambda) = 
 \prod_{j=1}^r (\cos \theta - \cos \theta_{r}) G_{\lambda}(\theta).$$
 On suppose de plus qu'il existe un unique r\'eel  $\theta_{0}$ dans $[0,2\pi[$ tel que 
   $f(\theta_{0})= \min_{\theta\in [0,2 \pi[} f(\theta)$. \\
   
  Alors pour $N$ assez grand les  valeurs propres de $T_{N}(f)$ sont de la forme $f\left(\frac{k\pi}{N+2}+
   \frac{\theta_{N,k}\pi}{N})\right)$ avec $\vert \theta_{N,k}\vert< 1$ .
 En cons\'equence 
  les plus petites valeurs propres de $T_{N}(f)$ sont de la forme $f\left(\frac{k\pi}{N+2}+
   \frac{\theta_{N,k}\pi}{N})\right)$ avec $\displaystyle{\lim_{N\rightarrow+ \infty}\frac{k\pi}{N+2} }= \theta_{0}$ et $\vert \theta_{N,k}\vert< 1$ .
   \end{theoreme}
  \begin{Preuve}{}
  On peut toujours se ramener \`a $f\ge0$.
  Si $m$ et $M$ d\'esignent respectivement le minimum et le maximum de $f$ sur le tore 
  $\mathbb T$, consid\'erons un r\'eel $\lambda$ v\'erifiant $m\le \lambda \le M$ 
  tel qu'aucun des ant\'ec\'edents de $\lambda$ par $f$ n'appartienne \`a 
  $J_{f} \cup \{0,\pi\}$. 
 % Avec les hypoth\`eses faites sur $P$ nous pouvons affirmer qu'il existe un polyn\^{o}me 
 % $Q$ de degr\'e $2 n_{0}$ tel que 
 %!TEX encoding = UTF-8 Unicode $\vert P\vert^2 -\lambda = \chi^{-n_{0}}\left( Q - \lambda \chi^n_{0}\right)$.
  Si $f_{1}$ est la fonction d\'efinie par la relation 
 $ f(\theta)= f_{1}(1 - \cos \theta)$,
   on notera 
  $\{\lambda'_{1}, \cdots, \lambda'_{r}\}$ l'ensemble $f_{1}^{-1} \{\lambda\}$. 
  Pour tout entier $j$, $1\le j \le r$ 
  on peut \'ecrire $ (1-\cos \theta)-\lambda'_{j}
= \frac{1}{2} \left( \vert 1-\chi \vert ^2- 2\lambda'_{j}\right)$. En posant 
$ \chi_{\lambda'_{j}} = (1-\lambda'_{j}) + i \sqrt {1-(\lambda_{j}'-1)^2}$ on peut \'ecrire 
cette \'equation sous la forme 
$ (1-\cos \theta)-\lambda'_{j}= -\frac{1}{2} \chi_{\lambda'_{j}}
(1-\overline{ \chi_{\lambda'_{j}}}\chi) (1-\overline{ \chi_{\lambda'_{j}}}\bar \chi).$
D'apr\`es les hypoth\`eses de notre \'enonc\'e, il existe  une fonction 
$H_{\lambda}(\theta)$ $2\pi$-p\'eriodique strictement positive sur 
$[0,2\pi[$, et $r$  complexes distincts tels que 
$$ f -\lambda =  \epsilon (-\frac{1}{2})^r\left( \prod _{j=1}^r 
  \chi_{\lambda'_{j}}
  (1-\overline{ \chi_{\lambda'_{j}}}\chi) (1-\overline{ \chi_{\lambda'_{j}}}\bar \chi)\right)
  H_{\lambda}$$
  avec $\epsilon \in \{-1,1\}$. 
   Nous pouvons finalement \'ecrire que
   $T_{N}(f-\lambda) =T_{N}(\tilde f_{\lambda})$
   avec 
$$
 \tilde f_{\lambda}= 
C_{1} (1-\overline{ \chi_{\lambda'_{1}}} \chi) (1- \overline{ \chi_{\lambda'_{1}}} \bar \chi) 
\cdots (1-\overline{ \chi_{\lambda'_{r}}} \chi) (1- \overline{ \chi_{\lambda'_{r}}} \bar \chi)
\frac{1}{\vert P_{N+r,\lambda}\vert ^2}
$$
o\`u $C_{1}$ est une constante et $P_{N+r,\lambda}$ le polyn\^{o}me pr\'edicteur 
de degr\'e $N+r$ assoc\'i\'e \`a la fonction $H_{\lambda}$.
 Pour $R\in ]0,1[$ posons maintenant 
 \begin{equation} \label{star}
f_{\lambda, R} = C_{1}\prod_{j=1}^r (1- R\overline{ \chi_{\lambda'_{j}}} \chi) (1-R \overline{ \chi_{\lambda'_{j}}} \bar \chi) 
 \frac{1}{\vert P_{N+r,\lambda} \vert ^2}.
\end{equation}
 Pour faire un calcul d'inverse pr\'ecis nous allons utiliser le corollaire \ref{COCO1} et le lemme  \ref{DEUX} qui est \'enonc\'ee ici  et d\'emontr\'ee dans l'appendice de cet article. 
 \begin{lemme}
 Si $f$ est une fonction r\'eguli\`ere appartenant \`a $\mathcal A (\mathbb T,s)$ 
 avec $s>\frac{3}{2}$, et $k$ un entier dans $[0,N]$, nous avons 
 $\beta_{k,N} = \overline{\widehat{g^{-1}(0)} }\widehat{g^{-1}(k)}+O(N^{-s})$,
 o\`u $g$ est la fonction de $\mathbb H^{+} $ telle que $f= g \bar g$, et o\`u
 $\beta_{k,N}$ est le coefficient de $\chi^k$ dans le lyn\^{o}me pr\'edicteur de 
 $f$.
 \end{lemme}
 D'apr\`es les hypoth\`eses il existe deux r\'eels 
 $0<\tau_{1}<1<\tau 2$ quel que soit $\lambda$ dans $[m,M]$ la fonction $H_{\lambda}$ est la restriction au tore d'une fonction strictement positive sur la couronne 
 $\{ z \vert \tau_{1}<\vert z \vert <\tau_{2}\}$ Le corollaire 
 \ref{COCO} permet donc d'\'ecrire qu'on peut consid\'erer que pour $N$ assez grand $P_{N+r, \lambda}$ est tr\`es proche d'un polyn\^{o}me de degr\'e fini $N_{0}$, l'entier $N_{0}$ \'etant ind\'ependant de $\lambda$ . En effet en utilisant la propri\'et\'e \ref{DEUX}
 on peut montrer que les $N_{0}$ premiers termes des polynmes pr\'edicteurs S'approchent de quantit\'es fixes quand 
 $N$ est assez grand. 
 Ces deux r\'esultats nous permettent d'\'ecrire que pour $\epsilon>0$ fix\'e il existe un entier $N_{\epsilon}$ que pour tout $\lambda'$ dans $[0,2]$ on ait un polyn\^{o}me de degr\'e 
 $N_{0}$, not\'e $P_{N_{0},\lambda}$ tel que 
$$\max_{\theta \in [0, 2 \pi[} \Bigl \vert  P_{N+r,\lambda} (e^{i\theta})
- P_{N_{0},\lambda} (e^{i\theta})
\Bigr \vert \le \epsilon.$$
Posons maintenant 
$$ f_{0,\lambda} = 
C_{1} (1-\overline{ \chi_{\lambda'_{1}}} \chi) (1- \overline{ \chi_{\lambda'_{1}}} \bar \chi) 
\cdots (1-\overline{ \chi_{\lambda'_{r}}} \chi) (1- \overline{ \chi_{\lambda'_{r}}} \bar \chi)
\frac{1}{\vert P_{N_{0},\lambda}\vert ^2}
$$
et 
$$ f_{0,\lambda,R} = 
C_{1} (1-R \overline{ \chi_{\lambda'_{1}}} \chi) (1- R\overline{ \chi_{\lambda'_{1}}} \bar \chi) 
\cdots (1-R \overline{ \chi_{\lambda'_{r}}} \chi) (1- R\overline{ \chi_{\lambda'_{r}}} \bar \chi)
\frac{1}{\vert P_{N_{0},\lambda}\vert ^2}
$$

Nous pouvons remarquer que si 
$T_{N} (f_{0,\lambda}) $ est inversible alors 
$T_{N} (\tilde f_{\lambda})$ l'est aussi.
En effet on a 
$\Vert T_{N} (\tilde f_{\lambda}) - T_{N} (f_{0,\lambda}) \Vert
\le 2 \pi \epsilon$ et 
$$ T_{N} (\tilde f_{\lambda}) = T_{N} (f_{0,\lambda})
\left( 1 -T_{N} (f_{0,\lambda})^{-1} \left( T_{N} (f_{0,\lambda}) -T_{N} (\tilde f_{\lambda}) \right) \right).$$

Nous allons maintenant pr\'eciser les coefficients  intervenant
dans le calcul de $T_{N}^{-1}\left( f_{0,\lambda}\right)_{k+1,l+1}$. Pour cela nous allons d'abord calculer les coefficients correspondants de $T_{N}^{-1}\left( f_{0,\lambda,R}\right)$.
Nous avons la d\'ecomposition 
$ f_{0,\lambda,R} = C_{1 }g_{1,0,\lambda}g_{2,0,\lambda}$ 
avec 
$g_{1,0,\lambda} = \displaystyle{\prod_{j=1}^r (1- R\overline{ \chi_{\lambda'_{j}}} \chi) \frac{1}{P_{N_{0},\lambda}}}$
et 
$g_{2,0,\lambda} = \displaystyle{\prod_{j=1}^r (1-R \overline{ \chi_{\lambda'_{j}}} \bar \chi) 
\frac{1}{\bar P_{N_{0},\lambda}}}$.
Le th\'eor\`eme \ref{UN} permet d'affirmer que 
$T_{N}\left( f_{0,\lambda}\right)$
est inversible en se souvenant que les racines de $P_{N_{0},\lambda}$ sont toutes de module strictement sup\'erieur \`a $1$ (voir \cite{Ld}).
En utilisant la d\'ecomposition en \'el\'ements simples des fractions rationnelles nous pouvons \'ecrire 
$$ \frac{1}{g_{1,0,\lambda}} = \sum_{j=1}^r \frac{A_{j}}{1-\omega_{j,\lambda,R}\chi}
 P_{N_{0},\lambda}, \quad
 \frac{1}{g_{2,0,\lambda} } = \sum_{j=1}^r \frac{A_{j}}{1-\omega_{j,\lambda,R}\bar \chi}
\bar P_{N_{0},\lambda},$$
 avec 
 $ \omega_{j,\lambda,R}= R \overline{ \chi_{\lambda'_{j}}}$. 
 Dans la suite, pour all\'eger les notations nous noterons par simplement 
 $\omega_{j}$ la quantit $ \omega_{j,\lambda,R}$, et par $P_{N_{0}}$ 
 le plolyn\^{o}me $P_{N_{0},\lambda}$.\\
 Posons maintenant
 $ \displaystyle{\pi_{+}\left( \frac{\chi^k}{ g_{2,0,\lambda}} \right) = \sum_{s=0}^k \gamma_{s} \chi^s},$
 la quantit\'e 
 $x_{k,R} = \pi_{+}\left( \pi_{+}\left( \frac{\chi^k }{g_{2,0,\lambda}}\right) \tilde \Phi_{N}\right)$
 s'\'ecrit 
 $ x_{k,R}= \displaystyle{\sum_{s=0} ^k \gamma_{s}
 \pi_{+} \left( \frac{\chi^{s-N-1}}{g_{1,0,\lambda}} 
 g_{2,0,\lambda}\right)}$.
  En calculant il vient 
$$ x_{k,R}= \sum_{s=0} ^k \gamma_{s} \sum_{j=1}^r
\pi_{+}\left(\frac{A_{j} \chi^{s-N-1} Q(\bar \chi)}
{(1-\omega_{j}\chi)}
\frac{P_{N_{0}}(\chi)}{\bar P_{N_{0}} (\bar \chi)}\right)$$
en posant 
$Q(z) = \displaystyle{\prod _{n=1}^r\left(1-\omega_{n}(z)\right)}$.
On obtient finalement : 
$$x_{k,R}= \sum_{j=1}^r A_{j}\left( \sum_{s=0} ^k \gamma_{s}
\omega_{j}^{N+1-s}\right)
\frac{Q(\omega_{j})}{(1-\omega_{j}\chi)}
\frac{P_{N_{0}}(\omega_{j}^{-1})}{\bar P_{N_{0}}(\omega_{j})}.$$
Si
 $x_{l,R} = \pi_{+}\left( \pi_{+}\left( \frac{\chi^l }{\bar g_{1,R}}\right) \bar \Phi_{N}\right)$, 
 en posant 
 $\displaystyle{\pi_{+}\left( \frac{\chi^l }{\bar g_{1,R}}\right)= \sum_{s=1} ^l \gamma'_{s} \chi^s}$
  un m\^{e}me calcul que pr\'ec\'edemment permet d'\'ecrire 
 $$ x_{l,R}= \sum_{j=1}^r \bar A_{j}\left( \sum_{s=0} ^l \gamma'_{s}
\omega_{j}^{N+1-s}\right)
\frac{Q(\bar \omega_{j})}{(1-\bar \omega_{j}\chi)}
\frac{\bar P_{N_{0}}(\omega_{j})}{ P_{N_{0}}(\omega_{j}^{-1})}.$$

 Notons maintenant, $H_{N,\lambda,R}$ l'op\'erateur 
 $H_{\tilde \Phi_{N,\lambda,R}}  H_{\Phi_{N,\lambda,R}}$, o\`u 
 $H_{\tilde \Phi_{N,\lambda,R}} $  et $H_{\Phi_{N,\lambda,R}}$ sont les deux op\'erateurs de Haenkel d\'efinis \`a partir de la d\'ecomposition (\ref{star}). Dans la suite on posera 
 $H_{N,\lambda}$ la limite quand $R$ tend vers $1$ par valeurs inf\'erieures de $H_{N,\lambda,R}$.
 Nous avons, en posant $ Q_{m}(z) =\displaystyle{ \prod_{n=1, n\neq m} ^{r} (1-\omega_{n} z)}$, pour tout entier $m$, 
 $1\le m \le r$, et $z \in \mathbb C$,
 \begin{align*}
 H_{\Phi_{N}} \left( \frac{1}{1-\omega_{j} \chi}\right)
 &= \pi_{-}\left( \frac{g_{1,0,\lambda}}{g_{2,0,\lambda}(1-\omega_{j}\chi)}  \chi^{N+1}
 \right)
 = \pi_{-}\left(\sum_{h=1}^r A_{h}
 \frac{Q_{j}(\chi) \chi^{N+1}}{(1-\omega_{h}\bar \chi)} 
 \frac{\bar P_{N_{0}}(\bar \chi)} {P_{N_{0}}(\chi)}\right)\\
 &=\sum_{h=1}^r A_{h} \frac{Q_{j}(\omega_{h}) \omega_{h}^{N+1}}{(1-\omega_{h}\bar \chi)} 
  \frac{\bar P_{N_{0}}(\omega_{h}^{-1})} {P_{N_{0}}(\omega_{h})}
\end{align*}
et 
\begin{align*}
H_{\tilde \Phi_{N}} \left( \frac{1}{1-\omega_{h} \bar \chi}\right)
 &= \pi_{+}\left( \frac{g_{2,0,\lambda}}{g_{1,0,\lambda}(1-\omega_{h}\bar \chi)} 
  \chi^{-N-1}\right)
 = \pi_{+}\left(\sum_{i=1}^r A_{i}
 \frac{Q_{h}(\bar \chi) \chi^{-N-1}}{(1-\omega_{i} \chi)} 
 \frac{ P_{N_{0}}( \chi)} {\bar P_{N_{0}}(\bar \chi)}\right)\\
 &=\sum_{h=1}^r A_{i} \frac{Q_{h}(\omega_{i}) \omega_{i}^{N+1}}{(1-\omega_{i} \chi)} 
  \frac{ P_{N_{0}}(\omega_{i}^{-1})} {\bar P_{N_{0}}(\omega_{i})}
\end{align*}
Si $\mathcal H_{N,\lambda,R}$ d\'esigne la matrice de 
$H_{N,\lambda,R}$ dans la base $\{\frac {1}{1-\omega_{1}\chi}, \cdots,\frac{1}{1-\omega_{r}\chi}\}$
les coefficients de $\mathcal H_{N,\lambda,R}$ sont donc
\begin{equation}\label{COEF}
\left( \mathcal H_{N,\lambda,R}\right)_{i,j} = 
A_{i} \sum_{h=1}^r A_{h} \omega_{h}^{N+2} \omega_{i}^{N+2} 
B_{h,j}(\omega_{h},\omega_{i}),\quad 1\le i,j\le r
\end{equation}
avec 
$$B_{h,j} (\omega_{h},\omega_{i})= Q_{j}(\omega_{h}) Q_{h}(\omega_{i}) 
\frac{\bar P_{N_{0}}(\omega_{h}^{-1}) P_{N_{0}} (\omega_{i}^{-1})}
{ P_{N_{0}}(\omega_{h}) \bar P_{N_{0}} (\omega_{i})}.$$
 Nous pouvons aussi \'ecrire :
  $ T_{N}^{-1}\left( f_{0,\lambda,R}\right)  T_{N}(f_{0,\lambda}) = T_{N}^{-1}\left( f_{0,\lambda,R}\right)  
  T_{N} \left( f_{\lambda,R}\right) + T_{N}^{-1}\left( f_{0,\lambda,R}\right)  
  \left(T_{N} (f_{0,\lambda})-f_{0,\lambda,R}) \right).$
Il est d'autre part clair que pour $N$ fix\'e, nous avons 
$\displaystyle{\lim_{R\to 1^-} T_{N}\left(f_{0,\lambda,}
-f_{0,\lambda,R} \right)}=0$. 
On en d\'eduit que si  $ \displaystyle{\lim_{R\rightarrow 1^-} }  T_{N}^{-1}\left( f_{0,\lambda,R}\right)$ existe, alors $\displaystyle{  \lim_{R\rightarrow 1^-} 
T_{N}^{-1}\left( f_{0,\lambda,R}\right) 
T_{N}\left(f_{0,\lambda})\right) = 1}$ et 
$\displaystyle{  \lim_{R\rightarrow 1^-} T_{N}\left(f_{0,\lambda}\right) 
T_{N}^{-1}\left( f_{0,\lambda,R}\right) 
= 1}$, c'est \`a dire que si cette limite existe on a
$ \displaystyle{\lim_{R\rightarrow 1^-} }  
T_{N}^{-1}\left( f_{0,\lambda,R}\right) =
 T_{N}^{-1}(f_{0,\lambda}).$
  Puisque les limites $\displaystyle{\lim_{R\rightarrow 1^{-}} x_{k,R}}$
 et $\displaystyle{\lim_{R\rightarrow 1^{-}} x_{l,R}}$ existent et sont finies, si 
  $(I-H_{N,\lambda})$ est inversible, alors  
  $\displaystyle{ \lim_{r\rightarrow 1^{-}} \langle \left (I - H_{N,\lambda,R} \right)^{-1} x_{k,R}\vert x_{l,R}\rangle}$ sera une somme finie de termes 
  $\frac{1}{1-\chi_{\lambda'_{j}}\chi_{\lambda'_{h}}}$, qui sont d\'efinis puisque  les hypoth\`eses sur $\lambda$ impliquent que 
 $\chi_{\lambda'_{j}}\chi_{\lambda'_{h}}\neq 1$ (remarquons que le fait que $f^{-1} (\lambda) \not\in \{0,\pi\}$ implique que 
 pour tout $i$, $1\le i \le r$ $\chi_{\lambda_{i}}\neq -1,1$)
  
 On conclut qu'un r\'eel $\lambda \notin f^{-1} (J)$ sera une valeur propre de $T_{N}(f_{0,\lambda})$
 si et seulement si le $\mathbb C$-uplet $(\chi_{\lambda'_{1}}, \cdots, \chi_{\lambda'_{r}})$ est solution de l'\'equation :
 \begin{equation}\label{equation1}
 \det\left (I - H_{N,\lambda} \right)=0.
  \end{equation}
    
    %les coefficients $(\mathcal H_{N,\lambda})_{i,j}$ de la matrice  
   % $\mathcal H_{N,\lambda}$ peuvent  s'\'ecrire sous la forme d'une somme de termes 
 % $(\bar\chi_{\lambda'_{i}})^{2N+4} A_{i,j}$et de termes 
  % $(\bar\chi_{\lambda'_{i}})^{N+2} A'_{i,j}$ o\`u les coefficients $A_{i,j}$ et 
 % $A'_{i,j}$ se calculent \`a partir du d\'eveloppement (\ref{COEF}). 
  En d\'eveloppant $D_{N,\lambda}=\det\left (I - H_{N,\lambda} \right)$ par rapport \`a la premi\`ere ligne, on obtient une expression de la forme 
  \begin{equation} \label{DET1}
  D_{N,\lambda}= 1- \sum_{i=0}^{r}(\bar\chi_{\lambda'_{i}})^{2N+4} A_{i,i}
  +d_{0}(\chi_{\lambda'_{1}},\cdots,\chi_{\lambda'_{r}})
  \end{equation}
  o\`u $d_{0}$ est une fonction en $\chi_{\lambda'_{1}},\cdots,\chi_{\lambda'_{r}}$ 
  et 
 $\displaystyle{A_{j,j}= A_{j}^{2} 
% \prod_{s =1,s\neq j}^{r} \left( \bar \chi_{\lambda'_{j}}-\bar \chi_{\lambda'_{s}}\right)} 
 Q_{j}^2(\omega_{j})
 \frac{\bar P_{N+r} (\chi_{\lambda'_{j}})P_{N+r} (\chi_{\lambda'_{j}})}
{\bar P_{N+r} (\bar \chi_{\lambda'_{j}})P_{N+r} (\bar \chi_{\lambda'_{j}})}}$, c'est \`a dire que les coefficients $A_{j,j}$ sont non nuls.
  L'\'equation (\ref{equation1}) se ram\`ene  donc \`a 
  \begin{equation} \label{equation2} 
  \bar \chi_{\lambda'_{1}}^{2N+4} = d (\chi_{\lambda'_{1}},\cdots,\chi_{\lambda'_{r}})
  \end{equation}
  o\`u $d$ est une fonction en $\chi_{\lambda'_{1}},\cdots,\chi_{\lambda'_{r}}$ et o\`u on a suppos\'e, pour fixer les id\'ees, que 
  $A_{1,1} \neq 0$.
  On a donc n\'ecessairement : 
  $$ \chi_{\lambda'_{1}} = e^{i( \frac{\theta_{N,\lambda}}{2N+4} + \frac{2k\pi}{2N+4})}$$
  avec  $0\le k \le 2N+3$ et o\`u $\theta_{N,\lambda}$ d\'esigne une d\'etermination de l'argument de $d(\chi_{\lambda'_{1}},\cdots,\chi_{\lambda'_{r}})$ 
  (on prend une d\'etermination de l'argument 
  comprise entre $[0,2\pi[$).  De plus l'\'ecriture 
\begin{equation}\label{JOJO}
 \chi_{\lambda'_{1}} = (1-\lambda'_{1}) + i \sqrt {1-(\lambda'_{1}-1)^2}
 \end{equation}
implique qu'en fait $0\le k \le N+1$.
On obtient donc finalement : 
\begin{equation}\label{TONTON}
 \lambda'_{1}= 1- \cos \left(\frac{\theta_{N,\lambda}}{2N+4} + \frac{k\pi}{N+2} \right),
 \quad 0\le k \le N+1
 \end{equation}
 ce qui conduit au r\'esultat. 
  \end{Preuve}
  
  \section{Une application \`a l'inversion de certaines matrices de Toeplitz bandes}
   \subsection{Enonc\'e du th\'eor\`eme}
  Dans cette partie on \'etudie des matrices de Toeplitz $T_{N}$ 
 telles qu'il existe un entier naturel $n_{0}$ v\'erifiant la condition :
 $\vert k-l\vert >n_{0} \Rightarrow \left(T_{N}\right)_{k,l} =0,$ et l'on suppose $a_{j}= \bar a_{-j}.$
 pour tout entier $j$.
 Pour tout entier $j$ compris entre $-n_{0}$ et $n_{0}$, on pose :
 $a_{j}= \left( T_{N}\right)_{k, k+j}$.
 La matrice $T_{N}$ admet alors pour symbole la fonction $\phi$ d\'efinie sur le tore par 
 $$ \phi (e^{i \theta}) = \sum_{n=-n_{0}}^{n_{0}} a_{n}e^{i n \theta}.$$
 On pose pour tout complexe $z$ : 
 $\displaystyle{K (z) =  \sum_{n=-n_{0}}^{n_{0}} a_{m} z^{(m+n_{0})}.}$
 Dans la suite, on supposera que  $K$ n'admet que des racines de module diff\'erent de $1$ et que $K$ admet $n_{0}$ racines de module inf\'erieur \`a $1$ et 
 $n_{0}$ racines de module sup\'erieur \`a 1. 
 On notera $\alpha_{i}$ les racines de $K$ qui sont de module strictement 
 inf\'erieur \`a 1, compt\'ees avec leur ordre de multiplicit\'e not\'e $s_{i}$. 
  Nous avons bien s\^{u}r 
 $ K(\chi) = C \displaystyle {\prod_{i=1}^{\sigma} (\chi-\alpha_{i})^{s_{i}} 
 \prod_{j=1}^{\sigma } (\chi-\overline{\alpha^{-1}_{j}})^{s_{j}}}$ avec 
 $\displaystyle{\sum_{i=1}^\sigma s_{i}}=n_{0}$. Dans la suite on supposera $C=1$ et
 on \'ecrira
 $ \phi(\chi) =  g_{1}g_{2}$ avec 
$$ g_{1} = (-1)^{j} \prod_{j=1}^{\sigma } \overline{\alpha^{-1}_{j}}
\prod_{h=1}^{s_{j}} (1-\bar \alpha_{j} \chi)^{m_{j}},
 \quad 
 g_{2}= \prod_{i=1}^{\sigma } (1-\alpha_{i} \bar \chi)^{s_{i}}.$$
 Nous avons : 
 $$ \frac{1}{g_{1}} = \sum_{j=1}^{\sigma} \left( \sum_{h=1}^{s_{j}} 
 \frac{K_{h,j}}{ (1-\alpha_{j}\chi)^h}\right),
 $$ 
  et 
 $$ \frac{1}{g_{2}} = \sum_{i=1}^{s} \left( \sum_{t=1}^{s_{i}} 
 \frac{H_{t,i}}{ (1-\alpha_{i}\bar \chi)^t}\right).
  $$ 
 Nous pouvons maintenant donner l'\'enonc\'e suivant :
   \begin{theoreme}\label{theoreme4}
  Si $T_{N}$ est une matrice de Toeplitz bande v\'erifiant 
  les hypot\`eses pr\'ec\'edentes, nous avons 
   alors \begin{itemize}
   \item [i)]
   $T_{N}$ est inversible 
   \item [ii)]
   Si $\vert k-l\vert =N x $, $0<x<1$  alors pour tout 
   r\'eel $a\in ]0,1[$ on a 
   $\left( T_{N}^{-1}\right)_{k+1,l+1} = o(\rho^{a \vert k-l\vert })$
  pour tout r\'eel $\rho$ tel que $1>\rho> \max \{\vert \alpha_{i}\vert, 1\le i \le s\}.$
         \end{itemize}
  \end{theoreme}
    Dans la suite de ce travail, pour tous les entiers $n_{i}$ nous d\'esignerons par $E_{n}$     le sous espace vectoriel des fonctions de $\mathbb T$ dans $\mathbb R$
     engendr\'e par les fractions rationnelles $\frac{1}{(1-\bar\alpha_{j} \chi)^n}$ pour $1\le n\le m_{j}$
     $1\le j\le m$. Nous pouvons maintenant \'enoncer et d\'emontrer un lemme 
    pr\'eparatoire.
     \subsection{Un lemme pr\'eparatoire}
    \begin{lemme} \label{lemme1}
    Avec les hypoth\`eses et les notations d\'efinies ci-dessus, on d\'efinit $E$ comme la somme des 
   $ E_{n_{i}}$, $1\le i\le s$. Alors $E$ est stable par 
    $H_{\tilde \Phi_{N}} H_{\Phi_{N}}$.
        \end{lemme}
  \begin{Preuve}{}
    Nous devons donc calculer 
    $H_{\tilde \Phi_{N}}H_{\Phi_{N}}(\psi)$ avec 
    $\psi = \frac{1}{(1-\bar\alpha_{j}\chi)^{n}}$, 
    $1\le n \le m_{j}$  un \'el\'ement de $E$.
    On peut se ramener \`a $\psi = \frac{1}{(1-\tilde \beta_{1}\chi)^{n}}$, 
    $1\le n \le m_{1}$.
    On a 
    $$
 H_{\Phi_{N}}  (\psi) =\pi_{-}\left(B
   \chi^{N+1+S} \frac{\left(1-\bar \alpha_{1}\chi)^{s_{1}-n} 
   \cdots(1\bar \alpha_{m}\chi)^{s_{\sigma}}\right) }
   { \left(1-\alpha_{1}\bar \chi)^{s_{1}} \cdots 
  (1-\alpha_{s}\bar \chi)^{s_{\sigma}} \right)}\right).
    $$
      Il est clair que, via la d\'ecomposition en 
    \'el\'ements simples de 
   $  
    \frac{1 }
   { \left(1-\alpha_{1}\bar  \chi)^{s_{1}} \cdots 
  (1-\alpha_{s}\bar \chi)^{s_{\sigma}} \right) },
   $ la quantit\'e 
     $H_{\Phi_{N}} (\psi)$ est 
    combinaison lin\'eaire de termes 
    $\pi_{-}\left( \frac{\chi^r}{(1-\alpha_{i}\bar \chi )^m}\right) $
      avec $1\le m\le s_{i}$, $1\le i\le \sigma$, $r=O(N)$.
         $H_{\Phi_{N}} (\psi)$ est combinaison lin\'eaire de termes 
       $\pi_{-}\left( \frac{\chi^r}{(1-\alpha_{i}\bar \chi )^m}\right)$.
       \\
      On peut remarquer que  
 $\frac{1}{(1-\bar \chi \alpha_{i})^m} = 
    \displaystyle { \sum_{u=0}^\infty   {\alpha}_{i} ^{u} {\bar \chi} ^u \tau_{m}(u)}$ avec $\tau_{m}(u) = (u+m-1)\cdots (u+1)$ si $m>1$ et
    $\tau_{m}(u)=1$ si $m=1$.\\
On obtient alors facilement que 
    $$\pi_{-}\left( \frac{\chi^r}{(1- \alpha_{i} \bar \chi )^m}\right) = 
    \sum_{u> r }   {\alpha}_{i} ^{u} {\bar \chi} ^{u -r} 
    \tau_{m}(u)  =
   \sum_{w=1}^\infty   {\alpha}_{i} ^{w+r} {\bar \chi} ^w 
   \tau_{m}(w+r).$$ On peut montrer    
   qu'il existe $m$ polyn\^{o}mes $\phi_{k,m}$, $1\le k\le m$ de degr\'e $m-k$ 
   tels que
   $ \displaystyle{ \tau_{m}(w+r)= \sum_{k=1}^{m} \phi_{k,m}(r) \tau_{k}(w)}$.
   Ceci donne
   $$ \pi_{-}\left( \frac{\chi^r}{(1- \alpha_{i} \bar \chi)^m}\right) =  \alpha_{i}^{r+1} \bar \chi
  \left(  \sum_{k=1}^{m} \phi_{k,m}(r) \sum_{w=0}^\infty \alpha_{i}^{w} \chi^w 
   \tau_{k}(w) \right)
   $$ ou encore, pour tout $i$, $1\le i\le s$,
   $$  \pi_{-}\left( \frac{\chi^r}{(1-\alpha_{i}\bar \chi )^m}\right) = 
\alpha_{i}^{r+1} \bar \chi   \sum_{k=1}^m \phi_{k,m}(r)  \frac{1}{(1-\alpha_{i}\bar \chi )^k}.
   $$
   De la m\^{e}me mani\`ere, on v\'erifie que pour tout $m\le m_{j}$, 
   $H_{\tilde \Phi_{N}} \left( \frac{1}{(1-\alpha_{i}\bar \chi )^m}\right)$
    est combinaison lin\'eaire de termes 
     $\pi_{+}\left( \frac{\bar \chi^{r_{1}}}{(1- \bar\alpha_{j}\chi )^n}\right) $ avec 
    $1\le n\le m_{j}$, $1\le j\le m$, et $r_{1}=O(N)$.
    Un calcul identique au pr\'ec\'edent donne alors 
    $$ \pi_{+}\left( \frac{\bar \chi^r}{(1- \bar \alpha_{j}\chi )^n}\right) =
    \sum_{k=1}^m \alpha_{i}^{r} \phi_{k,n}(r) 
    \frac{1}{(1-\bar \alpha_{j}\chi)^k}$$
    ce qui d\'emontre le lemme.\\
    \end{Preuve}
 \subsection{La preuve du th\'eor\`eme.}
  \begin{Preuve}{} Le point i) du th\'eor\`eme est une application imm\'ediate du th\'eor\`eme \ref{UN}. \\
  Pour le point ii) 
  posons  $\left( T_{N}^{-1}\right)_{k+1,l+1}=T_{1,k,l}+T_{2,k,l}$, les quantit\'es  
 $T_{1,k,l}$ et $T_{2,k,l}$ \'etant respecivement \'egales \`a 
   $T_{1,k,l}= \langle \pi_{+} \left( \frac{\chi^k}{g_{2}} \right) \vert
    \pi_{+}\left( \frac{\chi^l}{\bar g_{1}} \right)\rangle$ et, en reprenant les notations du lemme \ref{UN} \\
    $ T_{2,k,l}= \Bigl \langle\left(I-H_{\Phi_{N}} H_{\tilde \Phi_{N}}\right)^{-1}
    \left (\pi_{+}\left (\pi_{+} \left( \frac{\chi^k}{g_{2}} \right) \tilde \Phi_{N}\right)\right)
  \Big\vert 
    \pi_{+}\left(  \pi_{+}\left( \frac{\chi^l}{\bar g_{1} }\right) \overline{ \Phi_{N}}\right)
   \Bigr \rangle.$\\
    Les calculs pr\'ec\'edents,  via les d\'ecompositions en \'el\'ements simples de $\frac{1}{g_{1}}$ et $\frac{1}{g_{2}}$, 
    donnent que la quantit\'e 
     $T_{1,k,l}$ vaut
    $$ T_{1,k,l}=\sum_{1\le i,j\le \sigma} \left(\sum _{1\le t\le s_{i}, 1\le h\le s_{j}} 
   H_{t,i}  \bar K_{h,j} \Bigl\langle \pi_{+}\left(\frac{\chi^k}{(1-\alpha_{i}\bar \chi)^{t}}\right)
   \Big\vert 
    \pi_{+} \left( \frac{\chi^l}{(1-\alpha_{j}\bar \chi)^{h}}\right)\Bigr\rangle\right),$$ 
    ce qui est aussi, en r\'eutilisant les calculs effectu\'es dans le lemme \ref{lemme1} 
     $$ T_{1,k,l}=\sum_{1\le i\le s, 1\le j \le m} \left(\sum _{1\le t\le s_{i}, 1\le h\le m_{j}}
     H_{t,i}  \bar K_{t,j} S_{1,k,l,i,j}(t,h)\right),$$ 
    avec 
      $$S_{1,k,l,i,j} (t,h)=\Bigl\langle\displaystyle{\sum_{u=0}^k \tau_{t}(u)
     \alpha_{i}^{u}  \chi^{k-u} \Big\vert \sum_{v=0}^l \tau_{h} (v)
     { \alpha_{j}}^{v} \chi^{l-v}}\Bigr \rangle.$$
     On obtient, si on suppose $l\ge k$, 
          $$S_{1,k,l,i,j} (t,h)= \sum_{u=0}^k \tau_{t}(u) \tau_{h}( l-k+u) 
\alpha_{i}^{u}{\bar  \alpha_{j}}^{l-k+u}.$$
      C'est \`a dire que
      $\vert S_{1,k,l,i,j} (t,h)\vert = O\left( l^m \left(\max_{j\in \{1,\cdots, \sigma\}}\left( \vert \alpha_{j}\vert \right)\right)^{\vert l-k\vert }\right)$,
      ce qui donne, pour tout r\'eel $a\in ]0,1[$  
      $\vert T_{1,k,l}\vert = o\left( \rho^{a \vert l-k\vert} \right),$ avec la d\'efinition de $\rho$ et de $a$.\\
      Nous pouvons d'autre part \'ecrire 
      $$ T_{2,k,l} = \langle S_{2,k} \vert S_{2,l}\rangle+\Bigl \langle 
      \left( \left(I-H_{\Phi_{N}} H_{\tilde \Phi_{N}}\right)^{-1}-I\right)
    S_{2,k}
    \vert 
   S_{2,l}
   \Bigr \rangle,$$        
      avec 
      $$ S_{2,k}= \pi_{+}\left( \pi_{+} (\frac{ \chi^k}{g_{2}})\tilde\Phi_{N}\right) $$
      et 
      $$S_{2,l} =  \pi_{+}\left( \pi_{+} \left(\frac{ \chi^l}{\bar g_{1}}\right)\bar\Phi_{N}\right).$$
      Remarquons tout d'abord que, via la d\'ecomposition 
      en \'el\'ements simples de $\frac{1}{g_{2}}$
       la quantit\'e $\pi_{+} (\frac{ \chi^k}{ g_{2}})$ est une combinaison lin\'eaire finie de termes 
      $\pi_{+} \left( 
      \frac{\chi^{k}}{(1-\alpha_{i}\bar \chi)^{s_{1}}}\right)$
      avec $0\le s_{1} \le \max \{s_{i} \vert 1\le i\le s\} $.
      Ecrivons 
      $$ \pi_{+} \left( 
      \frac{\chi^{k}}{(1-\alpha_{i}\bar \chi)^{s_{1}}}\right)
      = \sum_{u=0}^{k}  \alpha_{i}^u \chi ^{k-u}
      \tau_{s_{1}}(u).$$
      Nous constatons alors que
      $ S_{2,k}$ est combinaison lin\'eaire finie de termes   
      $ \displaystyle{\sum_{u=0}^{k}  \alpha_{i}^u       \tau_{\sigma}(u) 
      \pi_{+}\left ( \chi^{k-u-N-1}) \frac{g_{2}}{g_{1}}\right)}.$
      On v\'erifie que 
     $ \pi_{+}\left ( \chi^{k-u-N-1}) \frac{g_{2}}{g_{1}}\right)$
     est une combinaison lin\'eaire finie de termes 
     $\pi_{+} \left( \frac{\chi^{k-u-N-1+M}}
     {(1-\bar \alpha_{j} \chi)^\mu}\right)$
     o\`u $M$ est un entier born\'e ind\'ependamment de $N$.
   En \'ecrivant
     $$ \pi_{+} \left( \frac{\chi^{k-u-N-1+M}}
     {(1-\bar \alpha_{j})^\mu} \right)= \sum_{v\ge N+1+u-k} 
   {\bar \alpha_{j}}^v \chi^v \tau_{\mu}(v)$$
 on obtient   
 que $S_{2,k}$ est une somme finie de termes 
 $${\bar \alpha_{j}}^{N+1-k} \sum_{u=0}^{k}  
 \alpha_{i}^u \tau_{\sigma}(u) \sum_{v>N+1+u+s-k}
   { \bar \alpha_{j}}^{v-N-1+k} \chi^v \tau_{\mu}(v).$$
      En estimant de m\^{e}me  
      $S_{2,l}$ on obtient qu'il existe un entier $H$, ind\'ependant des entiers $k,l,N$ tel que que le produit scalaire 
      $\langle S_{2,k}\vert S_{2,l}\rangle$ est une combinaison
      lin\'eaire finie de termes d'ordres $O( \rho^{2N-k-l})
      O(N^H) $, c'est \`a dire que ce produit scalaire est d'ordre $o( \rho^{\vert k-l\vert })$ si $k$ et $l$ v\'erifient les hypoth\`eses du th\'eor\`eme.

      Enfin, en reprenant la d\'emonstration du  lemme \ref{lemme1} on obtient que 
$ \Vert H_{N}\Vert = O(\rho^{2N})= o(\rho^{\vert l-k\vert })$, ce qui ach\`eve de prouver le point ii).
  \end{Preuve}
 \subsection{Application aux matrices de Toeplitz dont le symbole est une fonction r\'eguli\`ere}
On dira ici qu'une fonction r\'eguli\`ere continue $f$ v\'erifie l'hypoth\`ese $(\mathcal H)$ 
s'il existe deux r\'eels $\rho_{1}<1<\rho_{2}$ tels que $f$ soit la restriction \`a $\mathbb T$
d'une fonction continue strictement  positive sur $\{ z\vert \rho_{1}<\vert z \vert <\rho_{2}\}$.
 \begin{corollaire}\label {COCO1}
 Soit $f$ une fonction r\'eguli\`ere v\'erifiant l'hypoth\`ese $(\mathcal H)$. Si $\vert k-l\vert$ tend vers l'infini avec $N$, 
 alors $\left(T_{N}(f)\right)_{k,l}^{-1}= O\left(  \left( \frac{1}{\rho}\right)^{\vert l-k\vert}\right)$ pour tout $\rho \in ]1,\rho_{2}[$.
  \end{corollaire}
      \begin{Preuve}{}
    Appelons $m$ le minimum pour $\theta \in [0, 2\pi[$ de $ f(e^{i\theta})$. Donnons-nous un r\'eel $\epsilon$ strictement positif inf\'erieur \`a $\frac{m}{2}$
      Soit $P$ un polyn\^{o}me \`a coefficients complexes tel que 
      $\Bigl \vert \vert P(z)\vert ^2 - f(z)  \Bigr \vert \le\min \left( \epsilon, \Vert T_{N}^{-1} (f) \Vert_{2}\right), $ pour tout $z$, 
      $\rho_{1} <z<\rho_{2}$. 
      Il est clair que le polyn\^{o}me $P$ n'admet pas de racines sur
       le tore.
     Nous pouvons \'ecrire, 
     $ T_{N} (\vert P \vert ^2) = T_{N}(f) - T_{N} (f) + T_{N}(\vert P\vert ^2) $
     ce qui est aussi 
       $$ T_{N} (\vert P \vert ^2) = T_{N}(f) \left[ 1+ T^{-1}_{N}(f) 
      \left( T_{N} (\vert P\vert ^2) - T_{N}(f)\right) \right],$$
      ce qui implique, gr\^{a}ce 
      aux hypoth\`eses faites sur $f$ et $P$, 
      $$ T_{N}^{-1} (\vert P \vert ^2) = \left[ 1+ T^{-1}_{N}(f) 
      \left( T_{N} (\vert P \vert ^2) - T_{N}(f)\right)\right]^{-1}T_{N}^{-1} (f),$$
Ce qui conduit \`a la majoration,
$$ \Vert  T_{N}^{-1} (\vert P \vert ^2) \Vert _{2}\le \Vert T_{N}^{-1} (f)\Vert_{2}
+ \Vert T_{N}^{-1} (f)\Vert_{2} ^2 \Vert  T_{N} (\vert P\vert ^2) - T_{N}(f)\Vert_{2} \frac{1}{1-\epsilon}$$
qui donne 
 $$ \Vert  T_{N}^{-1} (\vert P \vert ^2) \Vert _{2} 
\Vert  T_{N} (\vert P\vert ^2) - T_{N}(f)\Vert_{2}\le \frac{2\epsilon }{1+\epsilon}.$$
c'est \`a dire qu'en choisissant $\epsilon$ assez petit on 
peut affirmer que lamatrice$T_{N}^{-1} (\vert P \vert ^2) \left(T_{N} (\vert P\vert ^2) - T_{N}(f)
\right)$ est inversible.
On peut alors recommencer les calculs pr\'ec\'edents en inversant les rles de 
$T_{N} (\vert P \vert ^2)$ et $T_{N}(f)$,
ce qui conduit \`a 
     $ T_{N} (f) = T_{N}(\vert P\vert ^2) + T_{N} (f) - T_{N}(\vert P\vert ^2) $
     ou encore 
     $$ T_{N} (f) = T_{N}(\vert P\vert ^2) \left[ 1+ T^{-1}_{N}(\vert P\vert ^2) 
      \left( T_{N} (f) - T_{N}(\vert P\vert ^2)\right) \right]$$
    et finalement 
      $$ T_{N}^{-1} (f) = \left[ 1+ T^{-1}_{N}(\vert P\vert ^2) 
      \left( T_{N} (f) - T_{N}(\vert P\vert ^2)\right)\right]^{-1}T_{N}^{-1} (\vert P\vert ^2),$$
    ce qui est enfin 
          $$ T_{N}^{-1} (f) = T_{N}^{-1} (\vert P\vert ^2) + T^{-1}_{N}(\vert P\vert ^2) 
      \left( T_{N} (f) - T_{N}(\vert P\vert^2)\right) \left[ 1+ T^{-1}_{N}(\vert P\vert ^2) 
      \left( T_{N} (f) - T_{N}(\vert P\vert ^2)\right) \right]^{-1}T_{N}^{-1} (\vert P\vert ^2)$$
      ce qui donne le r\'esultat avec le th\'eor\`eme \ref{theoreme4},
      en remarquant que l'hypoth\`ese $(\mathcal H)$ implique 
      que $T^{-1}_{N}(\vert P\vert^2) = \left( \frac{1}{\rho}\right)^{\vert l-k\vert}$,
      pour tout $\rho\in ]1, \rho_{2}[$.
\end{Preuve}
    \subsection{Peuve du lemme  \ref{DEUX}}
Avec les notations de l'article et la formule habituelle d'inversion des matrices de Toeplitz \`a symbole r\'egulier  (voir\cite{JMR01}) on obtient
pour $f$ une fonction r\'eguli\`ere  $f = g \bar g$,  
$g \in \mathbb H^{+} (\mathbb T)$.
$$ \left( T_{N}^{-1}\right)_{l+1,k+1} = \langle \pi_{+}\left( \frac{\chi^{l}}{\bar g }\right)\vert \frac{\chi^{k}}{\bar g } 
\rangle - \langle \sum_{s=0}^{+\infty} \left( H^{*}_{\Phi_{N}} H_{\Phi_{N}}\right)^{s} \pi_{+} \bar \Phi_{N}
\pi_{+} \left( \frac{\chi^{l}}{\bar g } \right )\vert  \pi_{+} \bar \Phi_{N}
\pi_{+} \left( \frac{\chi^{k}}{\bar g }\right)\rangle.$$
avec $\Phi_{N} = \frac{g}{\bar g} \chi ^{N+1}$. Pour $l=0$cette formule nous donne
$$ \left( T_{N}^{-1}\right)_{1,k+1} = \langle \pi_{+}\left( \frac{1}{\bar g }\right)\vert \frac{\chi^{k}}{\bar g } 
\rangle - \langle \sum_{s=0}^{+\infty} \left( H^{*}_{\Phi_{N}} H_{\Phi_{N}}\right)^{s} \pi_{+} \bar \Phi_{N}
\pi_{+} \left( \frac{1}{\bar g }\right )\vert  \pi_{+} \bar \Phi_{N}
\pi_{+} \left( \frac{\chi^{k}}{\bar g }\right)\rangle.$$
Dans la suite de la preuve nous utiliserons les notations suivantes : 
$$\frac{1}{g}= \sum_{u\ge 0} \beta_{u }\chi^u \quad 
\frac{g}{\bar g} = \sum_{u \in \mathbb Z} \gamma_{u }\chi^u .$$
Puisque l'hypoth\\`ese $f \in \mathcal A(\mathbb T,s)$ implique que  $\frac{1}{g}$ et
$\frac{g}{\bar g} \in  \mathcal A(\mathbb T,s)$ (see \cite{Ru} or \cite{Ka})nous avons deux constantes positives   $K$ et $K'$telle que 
$$ \vert \beta_{u}\vert \le \frac{K}{u^s} \quad \forall u\in \mathbb N^\star\quad 
\mathrm{et} \quad \vert \gamma_{u}\vert \le \frac{K'}{u^s} \quad \forall u\in \mathbb Z^\star.$$
On peut alors \'ecrire 
$$ \langle \pi_{+}\left( \frac{1}{\bar g }\right)\vert \frac{\chi^{k}}{\bar g } \rangle =
\bar \beta_{0} \beta_{k},$$
$$ \pi_{+} \bar \Phi_{N}
\pi_{+} \left( \frac{1}{\bar g } \right ) = \pi_{+} \left( \bar \Phi_{N} \bar \beta_{0}\right) = 
 \bar \beta_{0} \sum_{v\ge N+1 } \bar \gamma_{-v} \chi^{v-N-1},$$
 et  $$\Bigl \vert  \sum_{v\ge N+1 } \bar \gamma_{-v} \chi^{v-N-1}\Bigr \vert 
 =O(N^{-s}).$$
 De plus 
$$ \pi_{+} \bar \Phi_{N}
\pi_{+} \left( \frac{\chi^{k}}{\bar g }\right)= \sum_{w=0}^k \bar \beta_{w} \left( \sum_{v\ge N+1-k+w}
\bar \gamma_{v}\chi^{v-N-1+k-w}\right).$$
Et d'autre part si  $\psi = \displaystyle{ \sum_{w\ge 0} \alpha_{w}\chi^w }$est une fonction dans  $ \mathbb H^+$ nous avons, en utilisant la continit\'e de la projection  $\pi_{-}$, 
$$ H_{\Phi_{N}} (\psi) = \sum_{w\ge 0} \alpha_{w} \left( \sum_{v>N+1+w} \gamma_{-v} \chi^{v+w+N+1}\right)$$
ce qui implique
\begin{align*}
 \Vert H_{\Phi_{N}} (\psi) \Vert _{2} &\le \sum_{w\ge 0} \vert \alpha_{w} \left (\sum_{v>N+1+w } \vert \gamma_{-v} \vert 
\right)\\ 
& \le \Vert \psi \Vert _{2} \left(\sum_{w\ge 0} \left(\sum_{v>N+1+w } \vert \gamma_{-v} \vert\right)^2\right)^{1/2}\\
&\le K \Vert \psi \Vert _{2} (N+1)^{-s}
\end{align*}
ce qui signifie que $\Vert H_{\Phi_{N}} \Vert \le K (N+1)^{-s+3/2}$.
Nous obtenons alors facilement 
$\Vert H^\star_{\Phi_{N}} \Vert \le K' (N+1)^{-s+3/2}$ et nous pouvons \'ecrire
$$ \Vert \sum_{s=0}^{+\infty} \left( H^{*}_{\Phi_{N}} H_{\Phi_{N}}\right)^{s} \pi_{+} \bar \Phi_{N}
\pi_{+} \left( \frac{1}{\bar g }\right )\Vert_{2} \le 
\frac{K'} {\left(1- {K'}^2 (N+1)^{-s+3/2}\right)^2} ( N)^{-s+1}.$$
Puisque pour tout entier  $k$ fix\'e dans $[0,N]$  nous avons
$$\Vert \pi_{+} \bar \Phi_{N} \pi_{+}\left( \frac{\chi^k}{\bar g}\right)\Vert _{2} \le \sum_{w=0}^k \vert \beta_{w} \vert \left( \sum_{v\ge N+1-k+w}
\vert  \gamma_{v}\vert \right)=O(1)$$
et nous pouvons \'ecrire 
$$ \left( T_{N}^{-1}\right)_{1,k+1} = \bar \beta_{0} \beta_{k} + O \left((N)^{-s}\right)$$
tce qui est le r\'esultat attendu.

    \bibliography{Toeplitz4}

\end{document}